\documentstyle[12pt]{article}

\font\tenof=msym10 at 12pt
\def\R{\mbox{\tenof R}}
\def\Z{\mbox{\tenof Z}}

\def\case#1#2{{\textstyle{#1\over #2}}}
\def\cqfd{\qquad\qquad\vrule height 4pt depth 2pt width 5pt}

\newcommand{\A}{{\cal A}\bigl(G(N)\bigr)}
\newcommand{\one}{\mbox{\bf 1}}
\newcommand{\smallone}{\mbox{$\mathbf{\scriptstyle1}$}}
\newcommand{\ap}{a^{\dagger}}
\newcommand{\cR}{\mbox{$\mathcal{R}$}}
\newcommand{\cH}{\mathcal{H}}

\begin{document}
%
%
{\parindent=0cm
{\Large \bf Comment on ``Generalized $q$-oscillators and their Hopf
structures''}

\vspace{2cm}
{\large \bf C Quesne\footnote{Directeur de recherches FNRS}%
\footnote{e-mail: cquesne@ulb.ac.be} and N Vansteenkiste%
\footnote{e-mail: nvsteen@ulb.ac.be}}

\medskip
Physique Nucl\'eaire Th\'eorique et Physique Math\'ematique,
Universit\'e Libre de Bruxelles, Campus de la Plaine CP229, Boulevard du
Triomphe, B-1050 Brussels, Belgium

\vspace{3cm}
Short title: \textsl{Generalized $q$-oscillators}


\bigskip
Submitted to: \textsl{J. Phys. A: Math. Gen.}

\bigskip
Date: \today

\vspace{3cm}
{\bf Abstract.} In a recent paper (1994 {\sl J.\ Phys.\ A: Math.\ Gen.\ }{\bf
27}
5907), Oh and Singh determined a Hopf structure for a generalized
$q$-oscillator
algebra. We prove that under some general assumptions, the latter is, apart
from
some algebras isomorphic to su$_q$(2), su$_q$(1,1), or their undeformed
counterparts, the only generalized deformed oscillator algebra that supports a
Hopf
structure. We show in addition that the latter can be equipped with a universal
$\cR$-matrix, thereby making it into a quasitriangular Hopf algebra.\par}

\newpage
%
%
\noindent
In a recent paper (henceforth referred to as I and whose equations will be
quoted
by their number preceded by I), Oh and Singh~\cite{oh} studied the
relationships
among various forms of the $q$-oscillator algebra and considered the conditions
under which it supports a Hopf structure. They also presented a generalization
of
this algebra, together with its corresponding Hopf structure.\par
%
%
In the present comment, our purpose will be twofold. First, we plan to show
that
under some general assumptions about the coalgebra structure and the antipode
map, the generalized $q$-oscillator algebra considered by Oh and Singh is,
apart
from some algebras isomorphic to su$_q$(2), su$_q$(1,1), or their undeformed
counterparts, the only generalized deformed oscillator algebra~(GDOA) that
supports
a Hopf structure. Second, we shall provide the universal $\cR$-matrix for this
deformed algebra and prove that the corresponding Hopf algebra is
quasitriangular.\par
%
%
Let us introduce GDOA's as follows:
\par
\bigskip
\noindent{\it Definition.} Let $\A$ be the associative algebra generated by the
operators $\bigl\{\one,a,\ap,N\bigr\}$ and the function $G(N)$, satisfying the
commutation relations
\begin{equation}
   \left[N,\ap\right] = \ap \qquad \left[N,a\right] = - a \qquad
\left[a,\ap\right] =
   G(N)      \label{eq:algebra}
\end{equation}
and the Hermiticity conditions
\begin{equation}
   (a)^{\dagger} = \ap \qquad \left(\ap\right)^{\dagger} = a \qquad N^{\dagger}
= N
   \qquad \left(G(N)\right)^{\dagger} = G(N)          \label{eq:hermiticity}
\end{equation}
where $G(z)$ is assumed to be an analytic function, which does not vanish
identically.\par
\bigskip
For
\begin{equation}
   G(N) = [\alpha N + \beta + 1]_q - [\alpha N + \beta]_q
   = \frac{\cosh(\varepsilon(\alpha N + \beta + 1/2))}{\cosh(\varepsilon/2)}
   \label{eq:ohalgebra}
\end{equation}
where $\alpha$ and $\beta$ are some real parameters, and $q = \exp \varepsilon
\in \R^+$, $\A$ reduces to the generalization of the $q$-oscillator
algebra considered by Oh and Singh~\cite{oh}.\footnote[2]{Actually, Oh and
Singh
considered a slightly more general algebra, wherein the first two relations
in~(\ref{eq:algebra}) are also deformed by the introduction of a real
parameter~$\eta$. We shall not do so here, as this additional parameter can be
incorporated into the definition of~$N$ by renormalizing the latter.}\par
%
%
Note that the definition of $\A$ differs from the usual definition of
GDOA's~\cite{daska}, wherein both a commutation and an anticommutation
relations
\begin{equation}
   \left[a,\ap\right] = F(N+1) - F(N)  \qquad \left\{a,\ap\right\} = F(N+1) +
F(N)
   \label{eq:daska}
\end{equation}
are imposed in terms of some structure function $F(z)$, assumed to be an
analytic
function, positive on some interval $[0,a)$ (where $a \in \R^+$ may be finite
or
infinite), and such that $F(0) = 0$. As in I, the reason for considering only
the first
relation in~(\ref{eq:daska}) is that the two relations do not prove compatible
with a coalgebra structure.\par
%
%
As a consequence of its definition, the algebra~$\A$ has a Casimir operator
defined by $C = F(N) - \ap a$, where $F(N)$ is the solution of the difference
equation $F(N+1) - F(N) = G(N)$, such that $F(0) = 0$. The present definition
of
GDOA's is therefore equivalent to the usual one~\cite{daska} only in the
representation wherein $C=0$, i.e., in a Fock-type representation.\par
%
%
Let us now try to endow some of the algebras $\A$ with a coalgebra structure
and
an antipode map, making them into Hopf algebras $\cH$. For the coproduct,
counit and antipode, let us postulate the  following expressions :
\begin{eqnarray}
   \Delta \left(\ap\right) & = & \ap \otimes c_1(N) + c_2(N) \otimes \ap \qquad
       \Delta (a)  =  a \otimes c_3(N) + c_4(N) \otimes a  \label{eq:deltaa}
\\[0.1cm]
   \Delta (N) & = & c_5 N \otimes {\bf 1} + c_6 {\bf 1} \otimes N
       + \gamma {\bf 1} \otimes {\bf 1}  \label{eq:deltan} \\[0.1cm]
   \epsilon \left(\ap\right) & = & c_7 \qquad \epsilon (a)  =  c_8 \qquad
       \epsilon (N)  =  c_9  \label{eq:epsilon}  \\[0.1cm]
   S\left(\ap\right) & = & - c_{10}(N) \ap \qquad S(a) =  - c_{11}(N) a \qquad
       S(N)  =  - c_{12} N + c_{13} {\bf 1}  \label{eq:antipode}
\end{eqnarray}
where $ c_i(N)$, $i=1$, $\ldots$ , 4, 10,~11, are functions of $N$, and $c_i$,
$i=5$,
$\ldots$, 9, 12 ,~13, and $\gamma$ are constants to be  determined. Such
expressions generalize those found in I for $G(N)$ given
by~(\ref{eq:ohalgebra}),
which correspond to
\begin{eqnarray}
   c_1(N) & = & \left( c_2(N) \right)^{-1} = c_3(N) = \left( c_4(N)
\right)^{-1} =
      q^{\alpha (N + \gamma ) / 2 } \nonumber \\[0.1cm]
   c_5 & = & c_6 = c_{12} = 1 \qquad c_7 = c_8 = 0 \nonumber \\[0.1cm]
   c_9 & = & \case{1}{2} c_{13} = - \gamma  \qquad c_{10}(N) = \left( c_{11}(N)
       \right)^{-1} = q^{\alpha /2 } \nonumber  \\[0.1cm]
   \gamma & = & \frac {2 \beta + 1 }{2 \alpha } - {\rm i}\,
       \frac{(2 k + 1) \pi }{2 \alpha \varepsilon} \qquad k \in \Z.
\label{eq:ohhopf}
\end{eqnarray}
\par
%
%
To remain as general as possible, we shall not start by making any specific
assumption about $G(N)$, except that it satisfies eqs.~(\ref{eq:algebra})
and~(\ref{eq:hermiticity}). For the moment, we shall also disregard the
Hermiticity conditions~(\ref{eq:hermiticity}) and work with complex algebras.
Only at the end will conditions~(\ref{eq:hermiticity}) be imposed.\par
%
%
In order that equations (\ref{eq:deltaa})--(\ref{eq:antipode}) define a Hopf
structure, the so-far undetermined functions and parameters must be chosen in
such a way that $\Delta$, $\epsilon$, and~$S$ satisfy the coassociativity,
counit
and antipode axioms, given in~(I35), and that in addition, $\Delta$ and
$\epsilon$ be
algebra homomorphisms.\par
%
%
In accordance with eq.~(\ref{eq:ohhopf}), we shall start by assuming that in
eq.~(\ref{eq:deltan}), $\gamma$ takes a nonvanishing value. By substituting
eq.~(\ref{eq:deltan}) into the coassociativity axiom~(I35a), and taking into
account that $\Delta$ must be an algebra homomorphism, we directly obtain
\begin{equation}
   c_5 = c_6 =1.   \label{eq:c56}
\end{equation}
\par
%
%
To derive the corresponding conditions for $\ap$ and $a$, it is  useful to
expand the
functions $c_i(N)$, $i = 1$, $\ldots$,~4, of eq.~(\ref{eq:deltaa}) into power
series
\begin{equation}
   c_i(N) = \sum_{A=0}^{\infty} \frac{1}{A!}\, c_i^{(A)}(0) N^A
\label{eq:series}
\end{equation}
where $c_i^{(A)}(N)$ denotes the $A$th derivative of $c_i(N)$, and to apply the
relation
\begin{equation}
   \Delta c_i(N) = \sum_{A,B=0}^{\infty} \frac{1}{A! \, B!}\,
c_i^{(A+B)}(\gamma)
   N^A \otimes N^B \qquad \mbox{\rm if} \qquad \Delta(N) = N \otimes {\bf 1} +
   {\bf1} \otimes N  + \gamma {\bf 1} \otimes {\bf 1}.   \label{eq:deltac}
\end{equation}
We then obtain in a straightforward way that $c_i(N)$, $i = 1$, $\ldots$,~4,
must satisfy the equations
\begin{equation}
   c_i^{(A)}(0) c_i^{(B)}(0) = c_i^{(A+B)}(\gamma) \qquad i = 1, \ldots ,4
\qquad
   A,B = 0,1,2, \ldots .   \label{eq:cond1}
\end{equation}
\par
%
%
By substituting now eqs.~(\ref{eq:deltaa})--(\ref{eq:antipode}) into the counit
and
antipode axioms~(I35b) and~(I35c), we easily get
\begin{eqnarray}
   c_i(- \gamma ) & = & 1 \qquad i = 1, \ldots ,4   \label{eq:cond2}\\[0.1cm]
   c_7 & = & c_8 = 0 \qquad c_9 = - \gamma   \label{eq:c789}
\end{eqnarray}
and
\begin{eqnarray}
   c_1( - N +1-2 \gamma ) & = & c_2(N)\, c_{10}(N) \qquad c_2 ( - N  -2 \gamma
)  =
      c_1(N-1)\, c_{10}(N)   \label{eq:cond3}  \\[0.1cm]
   c_3 ( - N -1-2 \gamma ) & = & c_4(N)\, c_{11}(N) \qquad c_4 ( - N  -2 \gamma
)
      =  c_3(N +1)\, c_{11}(N)   \label{eq:cond4} \\[0.1cm]
   c_{12} & = & 1 \qquad c_{13} = - 2 \gamma   \label{eq:c1213}
\end{eqnarray}
respectively.\par
%
%
It remains to impose that the algebra and coalgebra structures are compatible.
By
applying $\Delta$ or $\epsilon$ to both sides of the first two equations
contained in
(\ref{eq:algebra}), we obtain identities, while by doing the same with the
third one
and using  equations similar to~(\ref{eq:series}) and~(\ref{eq:deltac}) for
$G(N)$,
we are led to  the conditions
\begin{equation}
   c_2(N +1) \otimes c_3(N)   =  c_2(N) \otimes c_3(N-1) \qquad
   c_4 (N) \otimes c_1(N+1)  =  c_4(N-1) \otimes c_1(N)   \label{eq:cond5}
\end{equation}
\begin{eqnarray}
   & & G^{(A-B)}(0) (c_1 c_3 )^{(B)}(0) + (c_2 c_4 )^{(A-B)}(0) G^{(B)}(0)
      = G^{(A )}(\gamma) \nonumber \\[0.1cm]
   & & \qquad A = 0,1,2, \ldots \qquad B= 0,1, \ldots A   \label{eq:cond6}
\end{eqnarray}
and
\begin{equation}
   G(- \gamma ) = 0 .   \label{eq:cond7}
\end{equation}
\par
%
%
We note that the Hopf axioms directly fix the values of all the constants
$c_i$, $i =
5$, $\ldots$, 9, 12,~13, in terms of the remaining one $\gamma$, but that the
seven functions $c_i(N)$, $i = 1$, $\ldots$, 4, 10,~11, and $G(N)$ are only
implicitly  determined by eqs.~(\ref{eq:cond1}), (\ref{eq:cond2}),
(\ref{eq:cond3}),
(\ref{eq:cond4}), (\ref{eq:cond5}), (\ref{eq:cond6}), and~(\ref{eq:cond7}). We
shall
now proceed to show that the latter can be solved to  provide explicit
expressions
for the yet unknown functions of $N$  in terms of $\gamma$ and of some
additional
parameters.\par
%
%
Considering first the two conditions in~(\ref{eq:cond5}), we immediately see
that they can only be satisfied if there exist some complex constants
$k_1$,~$k_2$, such that
\begin{eqnarray}
   c_1(N+1) & = & k_1 c_1(N) \qquad c_4(N ) = k_1^{-1} c_4(N-1)  \nonumber
      \\[0.1cm]
   c_2(N+1) & = & k_2 c_2(N) \qquad c_3(N ) = k_2^{-1} c_3(N-1).
      \label{eq:solcond5}
\end{eqnarray}
These relations in turn imply that
\begin{equation}
   c_1(N) = \alpha_1 \mbox{\rm e}^{\kappa_1 N } \qquad
   c_2(N) = \alpha_2 \mbox{\rm e}^{\kappa_2 N } \qquad
   c_3(N) = \alpha_3 \mbox{\rm e}^{-\kappa_2 N } \qquad
   c_4(N) = \alpha_4 \mbox{\rm e}^{-\kappa_1 N }   \label{eq:c1234}
\end{equation}
where $\kappa_1 = \ln k_1$, $\kappa_2 = \ln k_2$, and $\alpha_i$, $i = 1$,
$\ldots$,~4, are some complex parameters. The latter are determined by
condition~(\ref{eq:cond2}) as
\begin{equation}
   \alpha_1 = \mbox{\rm e}^{\kappa_1 \gamma} \qquad
   \alpha_2 = \mbox{\rm e}^{\kappa_2 \gamma} \qquad
   \alpha_3 = \mbox{\rm e}^{-\kappa_2 \gamma} \qquad
   \alpha_4 = \mbox{\rm e}^{-\kappa_1 \gamma}.   \label{eq:alpha}
\end{equation}
It is then straightforward to check that the functions $c_i(N)$, $i= 1$,
$\ldots$,~4,
defined by~(\ref{eq:c1234}) and~(\ref{eq:alpha}), automatically  satisfy
condition~(\ref{eq:cond1}).\par
%
%
By inserting now eqs.~(\ref{eq:c1234}) and~(\ref{eq:alpha}) into
conditions~(\ref{eq:cond3}) and~(\ref{eq:cond4}), we directly obtain the
following
explicit expressions for $c_{10}(N)$ and $c_{11}(N)$,
\begin{equation}
   c_{10}(N) = \mbox{\rm e}^{-(\kappa_1 + \kappa_2 )( N + \gamma ) + \kappa_1}
   \qquad c_{11}(N) = \mbox{\rm e}^{(\kappa_1 + \kappa_2 )( N + \gamma ) +
   \kappa_2 }.   \label{eq:c1011}
\end{equation}
\par
%
%
The same substitution performed in condition~(\ref{eq:cond6}) transforms the
latter into
\begin{eqnarray}
   & & (\kappa_1 - \kappa_2 )^B\, e^{(\kappa_1 - \kappa_2) \gamma} \,
      G^{(A-B)}(0) + (-1)^{A-B} (\kappa_1 - \kappa_2)^{A-B} \,
      e^{-(\kappa_1 - \kappa_2) \gamma}\, G^{(B)}(0) = G^{(A)}(\gamma)
\nonumber
      \\[0.1cm]
   & & \qquad A = 0,1,2, \ldots \qquad B= 0,1, \ldots ,A.   \label{eq:condG}
\end{eqnarray}
It can be easily shown by induction over $A$ that whenever $\kappa_1 \ne
\kappa_2$, the solution of recursion relation~(\ref{eq:condG}) is given by
\begin{equation}
   \begin{array}[b]{llll}
      G^{(A)}(0) &= &(\kappa_1 - \kappa_2)^A\, G(0) &\qquad\mbox{\rm if $A$ is
even}
         \\[0.3cm]
      &= & (\kappa_1 - \kappa_2)^A \coth \bigl((\kappa_1 - \kappa_2) \gamma
         \bigr) G(0) &\qquad\mbox{\rm if $A$ is odd}
    \end{array} \label{eq:GAzero}
\end{equation}
and
\begin{equation}
   \begin{array}[b]{llll}
      G^{(A)}(\gamma) & = &(\kappa_1 - \kappa_2)^A\, G(\gamma) &\qquad
          \mbox{\rm if $A$ is even} \\[0.3cm]
      & = & (\kappa_1 - \kappa_2)^A \coth \bigl(2(\kappa_1 - \kappa_2) \gamma
          \bigr) G(\gamma) &\qquad \mbox{\rm if $A$ is odd}
   \end{array}\label{eq:GAgamma}
\end{equation}
where
\begin{equation}
   G(\gamma) = 2 \cosh \bigl((\kappa_1 - \kappa_2) \gamma \bigr) G(0).
   \label{eq:Ggamma}
\end{equation}
{}From (\ref{eq:GAzero}) and the Taylor expansion of $G(N)$, we then obtain
\begin{equation}
   G(N) = G(0)\, \frac{\sinh \bigl((\kappa_1 - \kappa_2) (N + \gamma)
\bigr)}{\sinh
   \bigl((\kappa_1 - \kappa_2) \gamma \bigr)} \qquad \kappa_1 \ne \kappa_2 .
   \label{eq:GN}
\end{equation}
Such a function also satisfies (\ref{eq:GAgamma}) and~(\ref{eq:Ggamma}), as
well
as the remaining condition~(\ref{eq:cond7}). Equations
(\ref{eq:GAzero})--(\ref{eq:GN}) remain valid for $\kappa_1 = \kappa_2$
provided
appropriate limits are taken. In such a case, function~(\ref{eq:GN}) becomes
\begin{equation}
   G(N) = G(0) \left({\bf 1} + \frac{N}{\gamma }\right) \qquad \kappa_1 =
\kappa_2 .
   \label{eq:limitGN}
\end{equation}
\par
%
%
Had we taken $\gamma = 0$ instead of $\gamma \ne 0$ in~(\ref{eq:deltan}),
a similar analysis would have led to
\begin{equation}
   \begin{array}[b]{llll}
       G(N) & = & G^{(1)}(0) \,\frac{\displaystyle\sinh \bigl((\kappa_1 -
\kappa_2) N
           \bigr)}{\displaystyle\kappa_1 - \kappa_2} &\qquad \mbox{\rm if
$\kappa_1
           \ne \kappa_2$} \\[0.3cm]
       & = & G^{(1)}(0) \, N &\qquad \mbox{\rm if $\kappa_1 = \kappa_2$}
   \end{array}    \label{eq:GNzero}
\end{equation}
and a Hopf structure given by (\ref{eq:c56}), (\ref{eq:c789}),
(\ref{eq:c1213}),
(\ref{eq:c1234}), (\ref{eq:alpha}), and~(\ref{eq:c1011}), but where $\gamma$ is
set
equal to $0$. For an appropriate choice of $G^{(1)}(0)$ (obtained by
renormalizing
$\ap$ and~$a$ if necessary), such a form of
$G(N)$ corresponds to the  complex $q$-algebra sl$_q$(2) if $\kappa_1 \ne
\kappa_2$, and to sl(2) if $\kappa_1 = \kappa_2$~\cite{majid}.\par
%
%
The remaining step in the construction of algebras $\A$ with a Hopf structure
consists in imposing the Hermiticity conditions (\ref{eq:hermiticity}) on the
algebraic structure. They require that the function $G(N)$, defined
in~(\ref{eq:GN}),
(\ref{eq:limitGN}), or~(\ref{eq:GNzero}), be a real function of $N$. For the
latter
choice, we obtain the real forms of sl$_q$(2) or sl(2), namely su$_q$(2) and
su$_q$(1,1), or su(2) and su(1,1)~\cite{majid}. It remains to consider the
former
choices for $\gamma$ non real, since the real $\gamma$~case comes down to the
$\gamma = 0$ one by changing $N$ into $N + \gamma$. For such $\gamma$~values,
function~(\ref{eq:limitGN}) cannot be Hermitian. It therefore only remains to
consider the case where $G(N)$ is given by~(\ref{eq:GN}).\par
%
%
In such a case, the discussion of the hermiticity conditions is rather involved
as
$G(N)$ depends upon two complex parameters $\kappa_1 - \kappa_2$, and
$\gamma$, in  addition to the nonvanishing real parameter $G(0)$. By setting
\begin{equation}
   \kappa_1 = \xi_1 + \mathrm{i} \eta_1 \qquad
   \kappa_2 = \xi_2 + \mathrm{i} \eta_2 \qquad
   \kappa \equiv \kappa_1 - \kappa_2 = \xi + \mathrm{i} \eta \qquad
   \gamma = \gamma_1 + \mathrm{i} \gamma_2   \label{eq:xieta}
\end{equation}
where $\xi_1$, $\eta_1$, $\xi_2$, $\eta_2$, $\xi$, $\eta$, $\gamma_1$,
$\gamma_2 \in \R$, the function $G(N)$, defined in~(\ref{eq:GN}), can be
rewritten
as
\begin{equation}
   G(N) = G(0) \left(\alpha(N) + \mathrm{i} \beta(N) \right)
\label{eq:decompGN}
\end{equation}
where
\begin{eqnarray}
   \alpha(N) & = & \frac{a(N) c+b(N) d}{c^2+d^2} \qquad \beta(N)  =  \frac{b(N)
c-a(N)
       d}{c^2+d^2} \nonumber \\[0.1cm]
   a(N) & = & \sinh\left(A(N)\right) \cos\left(B(N)\right) \qquad b(N)  =
\cosh
       \left(A(N)\right) \sin\left(B(N)\right) \nonumber \\[0.1cm]
   c & = & \sinh C  \cos D \qquad d = \cosh C \sin D \nonumber \\[0.1cm]
   A(N) & = & \xi(N+ \gamma_1) - \eta \gamma_2 \qquad B(N)  =  \xi \gamma_2
       + \eta (N+\gamma_1) \nonumber \\[0.1cm]
   C & = & \xi \gamma_1 - \eta  \gamma_2 \qquad D = \xi \gamma_2 +\eta
       \gamma_1.   \label{eq:decomp}
\end{eqnarray}
Hence, $G(N)$ is a real function of $N$ if and only if
\begin{equation}
   \beta(N) = 0.   \label{eq:condhermite}
\end{equation}
Note that from the expressions of $\alpha(N)$ and $\beta(N)$ given in
(\ref{eq:decomp}), it is clear that the parameter values for which $c$ and $d$
simultaneously vanish should be discarded.\par
%
%
Condition~(\ref{eq:condhermite}) has now to be worked out by successively
combining the cases where $\gamma_1 = 0$ and $\gamma_2\ne 0$, or $\gamma_1
\ne 0 $ and $\gamma_2\neq 0$, with those where $\xi\ne 0$ and $\eta=0$,
$\xi = 0$ and $ \eta \ne 0$, or $\xi\ne 0$ and $\eta\ne 0$. For instance, if
$\gamma_1$, $\gamma_2$, $\xi \ne 0$, and $\eta=0$, equation
(\ref{eq:condhermite}) can be written as
\begin{equation}
   \cosh \bigl(\xi(N + \gamma_1)\bigr) \sin(\xi\gamma_2)
      \sinh(\xi\gamma_1) \cos(\xi\gamma_2 )
   = \sinh\bigl(\xi(N + \gamma_1)\bigr) \cos(\xi\gamma_2)
      \cosh(\xi\gamma_1) \sin(\xi\gamma_2).    \label{eq:case1}
\end{equation}
As both sides of this relation have a different dependence on $N$, they must
identically vanish. Since $\xi\ne 0$  by hypothesis, we must therefore have
either
$\sin(\xi\gamma_2) = 0$ or $\cos(\xi\gamma_2 ) =0$. The first condition leads
to
$\gamma_2 = k \pi / \xi$, $k \in \Z_0$, while the second one gives rise to
$\gamma_2 = (2 k + 1) \pi / (2 \xi)$, $ k \in \Z $.\par
%
%
Similarly, if we assume that $\gamma_1$, $\gamma_2$, $\eta \ne 0$, and $\xi=0$,
we obtain that equation (\ref{eq:condhermite}) is equivalent to
\begin{equation}
   \sin\bigl(\eta(N + \gamma_1)\bigr) \cos(\eta\gamma_1) =
   \cos\bigl(\eta(N + \gamma_1)\bigr) \sin(\eta\gamma_1)     \label{eq:case2}
\end{equation}
or, by using some trigonometric identities,
\begin{equation}
   \sin(\eta N) = 0.   \label{eq:case2a}
\end{equation}
As $\eta\ne 0$, this relation cannot be satisfied as an operator identity.\par
%
%
By proceeding in this way, one can easily show the following result:
\par
\bigskip
\noindent{\it Proposition 1.} The algebras~$\A$ that support a Hopf structure
of
type (\ref{eq:deltaa})--(\ref{eq:antipode}) and are not isomorphic to
su$_q$(2),
su$_q$(1,1), su(2), su(1,1), are determined by eqs.~(\ref{eq:c56}),
(\ref{eq:c789}),
(\ref{eq:c1213}), (\ref{eq:c1234}), (\ref{eq:alpha}), (\ref{eq:c1011}), and the
following conditions
\begin{eqnarray}
   G(N) &=& G(0)\, \frac{\cosh\bigl(\xi(N +
\gamma_1)\bigr)}{\cosh(\xi\gamma_1)}
       \qquad  G(0), \xi  \in \R_0 \qquad  \gamma_1 \in \R \nonumber \\[0.1cm]
   \kappa_1 - \kappa_2 &=& \xi \qquad \gamma = \gamma_1 + i\,
       \frac{(2k+1)\pi}{2\xi} \qquad  k \in \Z.   \label{eq:prop1}
\end{eqnarray}
\par
\bigskip
\noindent
{\it Remark.} The isomorphism referred to in the proposition is an algebra (not
a
Hopf algebra) isomorphism. One can indeed obtain algebras~$\A$ that have the
commutation relations and Hermiticity conditions of su$_q$(2), su$_q$(1,1),
su(2),
su(1,1), but more general expressions for the coproduct, the counit, and the
antipode.\par
%
%
Comparing the results of Proposition~1 with eqs.~(\ref{eq:ohalgebra})
and~(\ref{eq:ohhopf}), we notice that provided we set $\kappa_1 = - \kappa_2$,
the
Hopf algebra so obtained does coincide with that derived by Oh and Singh, the
relations between the two sets of parameters being given by
\begin{equation}
   G(0) = \frac{\cosh \left( \varepsilon (2 \beta + 1)/2 \right)}
   {\cosh (\varepsilon/2)} \qquad \xi = \alpha \varepsilon \qquad
   \gamma_1 = \frac{2\beta + 1}{2\alpha} .   \label{eq:equivalence}
\end{equation}
Hence, we have:
\par
\bigskip
\noindent{\it Corollary 2.} The only algebras~$\A$ that support a Hopf
structure of
type (\ref{eq:deltaa})--(\ref{eq:antipode}) and are not isomorphic to
su$_q$(2),
su$_q$(1,1), su(2), su(1,1), are isomorphic to those considered in~I.\par
\bigskip
\noindent{\it Remarks.} (1) The Hopf algebra obtained here is slightly more
general than that constructed by Oh and Singh, as it contains the additional
parameter $\kappa_1 + \kappa_2$.\par
\noindent (2) Some further generalizations of the coproduct given in
eqs.~(\ref{eq:deltaa}) and~(\ref{eq:deltan}), obtained by introducing
additional
functions of~$N$, fail to provide new Hopf algebras.\par
%
%
Let us now turn ourselves to the second point of this comment, namely the
construction of the universal $\cR$-matrix for the Oh and Singh Hopf
algebra.\par
%
%
By first omitting the Hermiticity conditions, one obtains the following result:
\par
\bigskip
\noindent{\it Lemma 3.} The complex Hopf algebras~$\cH$, defined by
eqs.~(\ref{eq:algebra}), (\ref{eq:deltaa})--(\ref{eq:antipode}),
(\ref{eq:c56}),
(\ref{eq:c789}), (\ref{eq:c1213}), (\ref{eq:c1234})--(\ref{eq:c1011}),
and~(\ref{eq:GN}), can be made into quasitriangular Hopf algebras by
considering
the element $\cR \in \cH \otimes \cH$, given by
\begin{eqnarray}
  \cR & = & X^{-2(N + \gamma \smallone) \otimes (N + \gamma \smallone)}
       \sum_{n=0}^{\infty} \frac{(1-X^2)^n}{[n]_X!} X^{-n(n-1)/2} Y^n
\lambda^{-2n}
       \nonumber \\
  & & \mbox{} \times \left((XY)^{(N + \gamma \smallone)} a\right)^n \otimes
       \left((XY)^{-(N + \gamma \smallone)} \ap\right)^n
\label{eq:complexR}
\end{eqnarray}
where
\begin{eqnarray}
  X & = & e^{(\kappa_1 - \kappa_2)/2} \qquad Y = e^{(\kappa_1 + \kappa_2)/2}
        \qquad \lambda^2 = - G(0) \frac{\sinh \left((\kappa_1 -
        \kappa_2)/2\right)}{\sinh \left((\kappa_1 - \kappa_2) \gamma\right)}
        \nonumber \\[0pt]
  [n]_X & = & \frac{X^n - X^{-n}}{X - X^{-1}} \qquad [n]_X! = [n]_X [n-1]_X
\ldots [1]_X
        \qquad [0]_X! = 1.     \label{eq:XY}
\end{eqnarray}
\bigskip
%
\noindent{\it Proof.} By direct substitution, one finds that $\cR$, defined
by~(\ref{eq:complexR}) and~(\ref{eq:XY}), satisfies the~relations
\begin{eqnarray}
  \left(\Delta \otimes \mbox{\rm id}\right) \cR & = & \cR_{13} \cR_{23} \qquad
          \left(\mbox{\rm id} \otimes \Delta\right) \cR = \cR_{13} \cR_{12}
          \nonumber \\[0.1cm]
  \tau \circ \Delta(h) & = & \cR \Delta(h) \cR^{-1}         \label{eq:qt}
\end{eqnarray}
where $\cR_{12}$, $\cR_{13}$, $\cR_{23} \in \cH \otimes \cH \otimes \cH$, and
for instance $\cR_{12} = \cR \otimes I$, while $\tau$ is the twist operator,
$\tau(a \otimes b) = b \otimes a$. \cqfd\par
%
%
By introducing now the additional conditions~(\ref{eq:prop1}) and $\kappa_1 +
\kappa_2 = 0$, and changing to Oh and Singh's notations~(\ref{eq:equivalence}),
we obtain the final result:
\par
\bigskip
\noindent{\it Proposition 4.} The Oh and Singh Hopf algebra, defined by
eqs.~(\ref{eq:algebra})--(\ref{eq:ohalgebra}),
(\ref{eq:deltaa})--(\ref{eq:ohhopf}),
is quasitriangular, with the $\cR$-matrix given by
\begin{eqnarray}
  \cR & = & q^{- \frac{1}{\alpha} \left[\left(\beta + \frac{1}{2}\right)^2 -
        \left(\frac{2k+1}{2\ln q}\pi\right)^2 + i \frac{(2\beta+1) (2k+1)}
        {2\ln q} \pi\right]} q^{-\alpha N \otimes N} \nonumber \\[0.1cm]
  & & \mbox{}\times \left(q^{- \left(\beta + \frac{1}{2} + i
        \frac{2k+1}{2\ln q} \pi\right) N} \otimes q^{- \left(\beta +
\frac{1}{2}
        + i \frac{2k+1}{2\ln q} \pi \right) N} \right) \nonumber \\[0.1cm]
  & & \mbox{}\times \sum_{n=0}^{\infty} \frac{\left[i\, (-1)^k \left(q^{1/2} +
        q^{-1/2}\right)\right]^n} {[n]_{q^{\alpha/2}}!} q^{-\alpha n(n-3)/4}
\left(
        \left(q^{\alpha N/2} a\right)^n \otimes \left(q^{-\alpha N/2}
\ap\right)^n
        \right).          \label{eq:realR}
\end{eqnarray}
\bigskip
\newpage
%
%

\end{document}